\documentclass[english,reqno]{amsart}
\usepackage[latin1]{inputenc}
\usepackage{graphicx,rotating}
\usepackage{amsthm,amsmath,amssymb}
\theoremstyle{plain}
\newtheorem{thm}{Theorem}[section]
\newtheorem{prop}[thm]{Proposition}
\newtheorem{lem}[thm]{Lemma}

\newtheorem{cor}[thm]{Corollary}
\theoremstyle{definition}
\newtheorem*{defn}{Definition}

\newtheorem{exa}[thm]{Example}
\theoremstyle{remark}
\newtheorem*{remark}{Remark}

\usepackage{graphics}

\newcommand{\eqdef}{\stackrel{\rm def}{=}}
\usepackage{babel}

\newcommand{\qc}{G(r,p,s,n)}

\newcommand{\inv}{\mathrm{inv}}

\newcommand{\sign}{\mathrm{sign}}
\newcommand{\Des}{\mathrm{Des}}
\newcommand{\HDes}{\mathrm{HDes}}

\newcommand{\maj}{\mathrm{maj}}
\newcommand{\Hilb}{\mathrm{Hilb}}

\newcommand{\fmaj}{\mathrm{fmaj}}
\newcommand{\fdes}{\mathrm{fdes}}
\newcommand{\ides}{\mathrm{ides}}
\newcommand{\ifmaj}{\mathrm{ifmaj}}
\newcommand{\icol}{\mathrm{icol}}

\newcommand{\s}{{\sigma}}

\newcommand{\Neg}{\mathrm{Neg}}
\newcommand{\des}{\mathrm{des}}
\newcommand{\col}{\mathrm{col}}

\def\s{ \sigma }

\def\ZZ{{\mathbb Z}}
\def\NN{{\mathbb N}}
\def\QQ{{\mathbb Q}}

\def\CC{{\mathbb C}}

\author{Riccardo Biagioli and Fabrizio Caselli}
\title[Enumerating  projective reflection groups]{Enumerating projective reflection groups}
\address{Universit\'e de Lyon, Universit\'e Lyon 1, Institut Camille
Jordan, UMR 5208 du CNRS, 69622 Villeurbanne, France\\
Dipartimento di matematica, Universit\`a di Bologna, Piazza di Porta San Donato 5, Bologna 40126, Italy}

\email{biagioli@math.univ-lyon1.fr, caselli@dm.unibo.it}

\begin{document}

\begin{abstract}
Projective reflection groups have been recently defined by the second author. They include a special class of groups denoted $G(r,p,s,n)$ which contains all classical Weyl groups and more generally all the complex reflection groups of type $G(r,p,n)$. In this paper we define some statistics analogous to descent number and major index over the projective reflection groups $G(r,p,s,n)$, and we compute several generating functions concerning these parameters. Some aspects of the representation theory of $G(r,p,s,n)$, as distribution of one-dimensional characters and computation of Hilbert series of invariant algebras, are also treated.
\end{abstract}

\maketitle


\section{Introduction}

The study of permutation statistics, and in particular of those depending on up-down or descents patterns, is a very classical subject of study in algebraic combinatorics that goes back to the early 20th century to the work of MacMahon \cite{mac60}, and has found a new interest in the last decade after the fundamental work of Adin and Roichman \cite{AR}. In their paper they defined a new statistic, the flag major index that extend to the Weyl group of type $B$ the concept of major index, which is a classical and well studied statistic over the symmetric group. This paper opened the way to several others concerning the study of statistics over classical Weyl groups, wreath products, and some types of complex reflection groups. Now we have a wide picture of several combinatorial and algebraic properties holding for all these families of groups. 

Recently, the second author introduced a new class of groups called projective reflection groups \cite{Ca1}. Among them, there is an infinite family denoted $G(r,p,s,n)$ where $r,p,s,n$ are integers such that $p$ and $s$ divide $r$, and $ps$ divides $rn$. They include the infinite family of irreducible complex reflection groups, in fact $G(r,p,1,n)=G(r,p,n)$. In this context is fundamental a notion of duality. Let $G=G(r,p,s,n)$, then his dual group is $G^*=G(r,s,p,n)$, where the parameters $p$ and $s$ are interchanged. In \cite{Ca1} it is shown that much of the representation theory of reflection groups can be extended to projective reflection groups, and that the combinatorics of the group $G=G(r,p,s,n)$ is strictly related to the invariant theory of the dual group $G^*$. 

In this paper we continue that study. We introduce two descent numbers $\des$ and $\fdes$, and a color sum $\col$ over $\qc$, that together with the flag major index $\fmaj$ \cite{Ca1} allow the extension to the groups $\qc$ of classical results in permutation enumeration. 

We show that the polynomial $\sum \chi(g) q^{\fmaj(g)}$, where $\chi$ denotes any one-dimensional character of $\qc$ and the sum is over $g \in \qc$, admits a nice product formula. This generalizes a classical result of Gessel and Simion \cite{W} for the symmetric group, and the main results of Adin-Gessel-Roichman \cite{AGR} for the Coxeter group of type $B$, and of the first author \cite{B1} for the type $D$ case. 
  
The enumeration of permutations by number of descents and major index yields a remarkable $q$-analogue of a well-known identity for the Eulerian polynomials, usually called Carlitz's identity. This identity has been generalized in several ways and to several groups, see e.g. \cite{ABR1, BB, BC, BZ, BZ1,CG, Re1}. All these extensions can be divided into two families, depending on the used descent statistic, either the number of geometric descents or the number of flag descents. We give a general method to compute the trivariate distribution of $\des$ (or $\fdes$), $\fmaj$ and $\col$ over $\qc$. This unifies and generalizes all previous cited results. 

We further exploit our method to compute the generating function of the six statistics $\des$, $\ides$, $\fmaj$, $\ifmaj$, $\col$, $\icol$ over the group $\qc$, ($\ides(g)$ denotes $\des(g^{-1})$ an similarly for the others). A specialization of this identity leads us to the computation of the generating function of the Hilbert series of the diagonal invariant algebras of the groups $\qc$. 

The definition of the above statistics depends on the particular order chosen on the colored integer numbers. By comparing our results with others in the literature, we deduce that a different choice of the order may lead to the same generating function. We conclude the paper by providing bijective explanations of this phenomenon.


\section{Notation and preliminaries}\label{notation}

In this section we collect the notations that are used in this paper as well as the preliminary results that are needed.

We let $\mathbb Z $ be the set of integer numbers and $\mathbb N$ be the set of nonnegative integer numbers. For $a,b\in \mathbb Z$, with $a\leq b$ we let $[a,b]=\{a,a+1,\ldots,b\}$ and, for $n\in \mathbb N$ we let $[n]\eqdef[1,n]$. For $r \in\mathbb N$, $r>0$ we let $R_r : \ZZ \rightarrow [0,r-1]$ the map ``residue module $r$", where $R_r(i)$ is determined by $R_r(i) \equiv i$ mod $r$. Moreover, we denote by $\zeta_r$ the primitive $r$-th root of unity $\zeta_r\eqdef e^{\frac{2\pi i}{r}}$. For $x \in \QQ$, we denote $\lfloor x \rfloor$ the biggest integer smaller than or equal to $x$. 
As usual for $n \in \NN$, we set $[n]_q\eqdef1+q+q^2+\ldots + q^{n-1}$, 
and ${\mathcal P}_n$ the set of nondecreasing sequences of nonnegative integers $(\lambda_1, \lambda_2, \ldots, \lambda_n)$, that is partitions of length less than or equal to $n$. 
Let $F \in \CC[[q_1,q_2,\ldots]]$ be a formal power series with complex coefficients, and $M$ a monomial. We denote by $\{F\}_{M}$ the power series obtained from $F$ discarding all its homogeneous components not divisible by $M$.
\smallskip

The main subject of this work are the complex reflection groups \cite{Sh}, or simply reflection groups. The most important example of a complex reflection group is the group of permutations of $[n]$, known as the symmetric group, that we denote by $S_n$.
We know by the work of Shephard-Todd \cite{ST} that all but a finite number of irreducible reflection groups are the groups $G(r,p,n)$, where $r,p,n$ are positive integers with $p|r$, that we are going to describe. If $A$ is a matrix with complex entries we denote by $|A|$ the real matrix whose entries are the absolute values of the entries of $A$.
The \emph{wreath product} $G(r,n)= G(r,1,n)$ is given by all $n\times n$ matrices satisfying the following conditions: 
the nonzero entries are $r$-th roots of unity; there is exactly one nonzero entry in every row and every column (i.e. $|A|$ is a permutation matrix). If $p$ divides $r$ then the {\em reflection group} $G(r,p,n)$ is the subgroup of $G(r,n)$ given by all matrices $A\in G(r,n)$ such that $\frac{\det A}{\det |A|}$ is a $\frac{r}{p}$-th root of unity. 

A \emph{projective reflection group}  is a quotient of a reflection group by a scalar subgroup (see \cite{Ca1}). Observe that a scalar subgroup of $G(r,n)$ is necessarily a cyclic group  $C_s$ of order $s$, generated by $\zeta_s I$, for some $s|r$, where $I$ denotes the identity matrix. It is also easy to characterize all possible scalar subgroups of the groups $G(r,p,n)$: in fact the scalar matrix $\zeta_sI$ belongs to $G(r,p,n)$ if and only if $s|r$ and $ps|rn$. In this case we let 
\begin{equation}\label{def-grpqn}
\qc\eqdef G(r,p,n)/ C_s.
\end{equation}
If $G=\qc$ then the projective reflection group $G^*\eqdef G(r,s,p,n)$, where the roles of the parameters $p$ and $s$ are interchanged, is always well-defined. We say that $G^*$ is the \emph {dual} of  $G$ and we refer the reader to \cite{Ca1} for the main properties of this duality. 
\smallskip

We sometimes think of an element $g\in G(r,n)$ as a \emph{colored permutation}. If the nonzero entry in the $i$-th row of $g\in G(r,n)$ is $\zeta_r^{c_i}$ we let $c_i(g)\eqdef R_r(c_i)$ and say that $c_1(g),\ldots,c_n(g)$ are the \emph{colors} of $g$. Now it is easy to see that an element $g\in G(r,n)$ is uniquely determined by the permutation $|g|\in S_n$ defined by  $|g|(i)=j$ if $g_{i,j}\neq 0$, and by its colors $c_i(g)$ for all $i \in[n]$. More precisely, when we consider an element of $G(r,n)$ as a colored permutation we represent its elements either as couples $(c_1,\ldots,c_n;\sigma)$, where $\s=\s(1)\cdots \s(n)$ is a permutation in $S_n$ and $(c_1,\ldots,c_n)$ is the sequence of its colors, or in {\em window notation} as
$$g=[g(1),\ldots, g(n)]=[\s(1)^{c_1},\ldots, \s(n)^{c_n}].$$
We call  $\sigma(i)$ the {\em absolute value} of $g(i)$, denoted $|g(i)|$.
When it is not clear from the context, we will denote $c_i$ by
$c_i(g)$. Moreover, if $c_i=0$, it will be omitted in the window notation of
$g$. In this notation the multiplication in $G(r,n)=\{(c_1,\ldots,c_n;\sigma) \mid c_i \in [0,r-1],\sigma \in
S_n\}$ is given by
$$(d_1,\ldots,d_n;\tau)(c_1,\ldots,c_n; \sigma)= (R_r(c_1+d_{\sigma(1)}),\ldots, R_r(c_n+d_{\sigma(n)});\tau\sigma).$$
If we let $\col(g)\eqdef c_1(g)+\cdots + c_n(g)$, then for $p|r$ we can simply express the reflection group $G(r,p,n)$ by
\begin{equation*}\label{def-grpn}
G(r,p,n)=\{(c_1,\ldots,c_n;\sigma)\in G(r,n) \mid \col(g)\equiv 0 \mod p\},
\end{equation*}
and $C_s$ as the cyclic group generated by $(r/s,\ldots,r/s; id)$.


\section{Descent-type statistics}\label{stat}

In this section we recall the notions of descents, flag descents, and flag major index for the wreath products. Then we introduce analogous definitions for descents and flag descents for projective reflection groups. 

\begin{defn} In all the paper we  use the following order, called {\em color order}
\begin{equation}\label{color-order}
1^{r-1}< \ldots < n^{r-1}< \ldots < 1^{1}< \ldots < n^1 < 0 < 1 < \ldots < n. 
\end{equation}
\end{defn}

For  $g \in G(r,n)$ we define the {\em descent set} as 
\begin{eqnarray}\label{des}
\Des_G(g)\eqdef \{i \in [0,n-1] \mid g(i) > g(i+1)\},
\end{eqnarray}
where $g(0):=0$, and denote its cardinality by $\des_G(g)$. A geometric interpretation of this set in terms of Coxeter-like generators and length can be given as in \cite[\S 7]{Re}. 
\smallskip
If we consider only positive descents we obtain 
\begin{equation}\label{desA}
\Des_A(g)\eqdef\Des_G(g) \setminus \{0\},
\end{equation}
and we denote $\des_A(g)$ its cardinality. The {\em flag major index} \cite{AR} and the {\em flag descent number} \cite{ABR1}, \cite{BB} are defined by
\begin{eqnarray}
\fmaj(g)&\eqdef& r \ \maj(g) +\col(g) \label{fmaj-wreath}\\
\fdes(g)&\eqdef& r \ \des_A(g) + c_1(g), \label{fdes-wreath}
\end{eqnarray}
where as usual the {\em major index} $\maj(g)=\sum_{i \in \Des_A(g)} i $ is the sum of all positive descents of $g$.
\smallskip

Following \cite[\S 5]{Ca1}, for $g=(c_1,\ldots,c_n;\sigma)\in \qc$ we let
\begin{eqnarray}
\HDes(g)&\eqdef&\{i\in [n-1] \mid c_i=c_{i+1} \textrm{, and }\sigma(i)>\sigma(i+1)\} \nonumber \\
h_i(g)&\eqdef&\#\{j\geq i \mid j\in \HDes(g)\} \label{def-pezzi}\\
k_i(g)&\eqdef& \left\{\begin{array}{ll} R_{r/s}(c_n)& \textrm{if }i=n\\k_{i+1}+R_r(c_i-c_{i+1})& \textrm{if }i\in [n-1].\nonumber
\end{array}\right.
\end{eqnarray}
We call the elements in  $\HDes(g)$ the \emph{homogeneous descents} of $g$. Note that the sequence $(k_1(g),\ldots,k_n(g))$ is a partition such that $g=(R_r(k_1(g)),\ldots,R_r(k_n(g)); \sigma)$. Moreover it is characterized by the following property of minimality: if  $\beta_1\geq \ldots \geq \beta_n$ and $g=(R_r(\beta_1),\ldots,R_r(\beta_n);\sigma)$, then $\beta_i\geq k_i(g)$, for all $i \in [n]$. 
\smallskip

For $g \in \qc$, we let
\begin{equation}
\lambda_i(g)\eqdef r\cdot h_i(g)+k_i(g).
\end{equation}
The sequence $\lambda(g)\eqdef(\lambda_1(g),\ldots,\lambda_n(g))$ is a partition. 
The {\em flag-major index} for $g \in \qc$ is defined by \cite[\S 5]{Ca1} 
\begin{equation}\label{def-fmaj-proj}
\fmaj(g)\eqdef |\lambda(g)|=\sum_{i=0}^n \lambda_i(g).
\end{equation} 
We define the {\em descent number} and the {\em flag descent number} of $g \in \qc$  respectively by
\begin{align}
\des(g)&\eqdef \Big\lfloor \frac{s \lambda_1(g)+r-s}{r} \Big\rfloor \quad {\rm and} \label{def-des-proj} \\
\fdes(g)&\eqdef \lambda_1(g) \label{def-fdes}.
\end{align} 
Finally, for $g=(c_1,\ldots,c_n;\sigma) \in \qc$ we define the {\em color} of $g$ by
\begin{equation}
\col(g) \eqdef \sum_{i=1}^n R_{r/s}(c_i).
\end{equation}
 
Note that the previous definitions do not depend on the particular representative of $g\in \qc$ chosen among  its $s$ lifts in $G(r,n)$.  
\begin{remark}
From \cite[Lemma 5.1]{Ca1}, it follows that for $p=s=1$ the flag major index on $\qc$ defined in \eqref{def-fmaj-proj} coincides with the flag major index of Adin and Roichman for wreath products $G(r,n)$ defined in \eqref{fmaj-wreath}. Moreover the definition of $\fdes$ is \eqref{def-fdes} is consistent with that in \eqref{fdes-wreath}, and so  
$$
\des(g)=\Big\lfloor\frac{\lambda_1(g)+r-1}{r}\Big\rfloor=\des_A(g)+\Big\lfloor\frac{c_1(g)+r-1}{r}\Big\rfloor
=\des_G(g),$$
where the last equality holds since
$$
\Big\lfloor\frac{c_1(g)+r-1}{r}\Big\rfloor=\left\{\begin{array}{ll} 0,&\textrm{if }c_1(g) = 0;\\1,&\textrm{otherwise}.\end{array}\right.
$$
These equalities give motivations to the definitions of the previous statistics.
\end{remark}
We finish this section with an example.
\begin{exa} Let $g=[2^2,7^3,6^3,4^5,8^1,1^1,5^3,3^2]\in G(6,2,3,8)$. We have $\HDes(g)=\{2,5\}$, $(h_1,\ldots,h_8)=(2,2,1,1,1,0,0,0)$ and $(k_1,\ldots,k_8)=(18,13,13,9,5,5,1,0)$, and so $\fdes(g)=6\cdot2+18=30$, $\des(g)=15$, and $\col(g)=6$.
\end{exa}
 

\section{One-dimensional characters and flag major index}\label{odcafmi}

The generating function of the major index with the unique nontrivial one-dimensional character  admits a nice factorization formula over the symmetric group, as shown in Theorem~\ref{gessim} by Gessel-Simion (see \cite{W}). The same happens for the other classical Weyl groups, as proved by Adin-Gessel-Roichman \cite{AGR} for the type $B$ case, and by the first author \cite{B1} for the type $D$ case. In this section we generalize these results to all projective groups of type $\qc$. 
\smallskip

The irreducible representations of $\qc$ are classified in \cite[\S 6]{Ca1}. In particular, for $n>2$, the one-dimensional characters of $\qc$ are all of the form 
$$\chi_{\epsilon, k}(g)= \epsilon^{\inv(|g|)} \zeta_r^{k \cdot c(g)},$$
where $\epsilon=\pm 1$, and $k\in [0, \frac{r}{p}-1]$ with the further condition that $s$ divides $kn$, and $c(g)$ is the sum of the colors of any element in $G(r,p,n)$ representing the class of $g$ (in particular $c(g)=\col(g)$ if $s=1$). As usual for $\s \in S_n$, we denote by $\inv(\s)\eqdef |\{(i,j) \in [n] \times [n] \mid i<j \ {\rm and} \ \s(i)>\s(j)\}|$, the number of its {\em inversions}, and by $\sign(\sigma)\eqdef(-1)^{\inv(\sigma)}$ its {\em sign}. Our main results is the following one.

\begin{thm}\label{maintwist}
Let $\chi_{\epsilon,k}$ be a one dimensional character of $\qc$. Then
$$
\sum_{g\in \qc} \chi_{\epsilon, k}(g)q^{\fmaj(g)}= \left[\frac{r}{p}\right]_{(\zeta^{k}q)^p}\left[\frac{2r}{p}\right]_{(\epsilon \zeta^{k}q)^p}\cdots \left[\frac{(n-1)r}{p}\right]_{(\epsilon^{n-2}\zeta^{k}q)^p}  \left[\frac{nr}{ps}\right]_{(\epsilon^{n-1} \zeta^k q)^{p}}\left\{[p]_{\zeta^k q}^{n-m}[p]_{\epsilon \zeta^k q}^{m} \right\}_{q^p}
$$
where
$m=\lfloor n/2 \rfloor$. 
\end{thm}

We recall that $ \{F\}_{q^p}$ is the polynomial obtained from $F$ discarding all the homogeneous components in the variable $q$ of degree not divisible by $p$. Hence for $r=2$ and $p=s=1$ in Theorem~\ref{maintwist}, we obtain \cite[Theorems 5.1, 6.1, 6.2]{AGR};  for  $r=s=2$ and $p=1$ \cite[Theorem 4.8]{B1}. 

In order to prove Theorem \ref{maintwist} we proceed in several steps passing through $G(r,n)$ and $G(r,1,s,n)$ until the general case of $\qc$. The basic stone is given by the following results for the symmetric group.
\begin{thm}[Gessel-Simion]\label{gessim} We have
$$
\sum_{\sigma \in S_n}\sign(\sigma)q^{\maj (\sigma)}=[1]_q[2]_{-q}[3]_q\cdots [n]_{(-1)^{n-1}q}.
$$
\end{thm}
The particular case of the Weyl groups  $B_n\eqdef G(2,n)$ is treated in \cite{AGR}. In this paper the authors focus on the sum 
$$
\sum_{g\in B_n}(-1)^{\ell(g)}q^{\fmaj(g)},
$$
where $\ell(g)$ is the Coxeter length of $g$ with respect to a given set of simple reflections, making use of a combinatorial interpretation of $\ell(g)$ due to Brenti \cite{bb05}. Theorem \ref{maintwist} is then achieved in this case since $(-1)^{\ell(g)}=(-1)^{\inv(|g|)}(-1)^{c(g)}$ for all $g\in B_n$. 
Although one can define an analogue of the Coxeter length for the wreath products $G(r,n)$ (see e.g., \cite{Ba},\cite{Re}), this does not lead to a one-dimensional character of the group, and the corresponding sum does not factorize nicely. This is why we focus on one-dimensional characters, obtaining in particular also a new proof in the case of Weyl groups of type $B$ that does not make use of the combinatorial interpretation of the length function.

The following is a technical lemma which is the multivariate extension of a result appearing in \cite[\S 5]{AGR}.
For $r>0$ and $k_{r-1},k_{r-2},\ldots,k_0\in \mathbb N$ such that $k_{r-1}+k_{r-2}+\cdots+k_0=n$ (i.e. $(k_{r-1},\ldots,k_0)$ is a composition of $n$) we let
$$
S(k_{r-1},\ldots,k_0)\eqdef\{\sigma\in S_n \mid \Des_A(\sigma)\subseteq\{k_{r-1}, k_{r-1}+k_{r-2},\ldots, k_{r-1}+k_{r-2}+\cdots+k_1\}\}.
$$
\begin{lem}\label{bin} We have
$$
\sum_{\sigma\in S(k_{r-1},\ldots,k_0)}\sign(\sigma)=\left\{ \begin{array}{ll}0, &\textrm{if at least 2 of the $k_i$'s are odd;}\\
\left( \begin{array}{c}
\lfloor \frac{n}{2}\rfloor\\ \lfloor\frac{k_{r-1}}{2}\rfloor, \ldots ,\lfloor \frac{k_{0}}{2}\rfloor \end{array}\right),& otherwise.
\end{array}
\right.
$$
\end{lem}

\begin{proof}
The window notation of $\sigma\in S(k_{r-1},\ldots,k_0)$ can be split into $r$ blocks of length $k_{r-1},\ldots,k_0$, the entries in each block being increasing. Consider the following involution on $S(k_{r-1},\ldots,k_0)$
$$
\left\{\begin{array}{ll}\sigma \mapsto (2i-1,2i)\sigma & \textrm{if $i$ is the smallest integer such that $2i-1$ and $2i$ are in different blocks,}\\ \sigma \mapsto \sigma & \textrm{if such $i$ does not exist.}
\end{array}
 \right.
$$
We recall that the window notation of $(2i-1,2i)\sigma$ is obtained from that of $\sigma$ by exchanging $2i-1$ and $2i$, and so this map is actually an involution on $S(k_{r-1},\ldots,k_0)$. Since  $\sigma$ and $(2i-1,2i)\sigma$ have opposite signs we can restrict our attention to those $\sigma\in S(k_{r-1},\ldots,k_0)$ such that $2i-1$ and $2i$ belong to the same block, for all $i\in \{1,\ldots,\lfloor \frac{n}{2} \rfloor\}$. It is clear that if two among the $k_i$'s are odd such $\sigma$ does not exist and the result follows in this case. 
We can thus assume that at most one of the $k_i$'s is odd, and in this case all entries of the window notation of $\sigma$ are uniquely determined by the even entries. This is because an odd entry $j$ different from $n$ belongs to the same block of $j+1$ and, if $n$ is odd, $n$ necessarily belongs to the unique block of odd length. The number of even entries is $\lfloor\frac{n}{2}\rfloor$ and these should be combined into $r$ parts of cardinality $\lfloor\frac{k_{r-1}}{2}\rfloor,\ldots, \lfloor \frac{k_{0}}{2}\rfloor$. It remains to show that all these elements have even sign. We first note that there are no inversions of the form $(\s^{-1}(2j-1),\s^{-1}(2j))$. Then we observe that, for all $i<2j-1<n$, $(\s^{-1}(2j-1),\s^{-1}(i))$ is an inversion if and only if $(\s^{-1}(2j),\s^{-1}(i))$ is itself an inversion. This completes the proof if $n$ is even. If $n$ is odd we only have to make the final remark that there are an even number of inversions of the form $(\s^{-1}(n),\s^{-1}(i))$, since $n$ appears in an odd position.
\end{proof}

\begin{thm}\label{grn}We have
$$
\sum_{g\in G(r,n)}\sign(|g|)q^{\fmaj(g)}=[r]_q[2r]_{-q}[3r]_q\cdots [nr]_{(-1)^{n-1}q}.
$$\end{thm}
\begin{proof}
We let $U \subset G(r,n)$ be given by $U=\{g\in G(r,n) \mid g(1)<g(2)<\cdots<g(n)\}=\{g\in G(r,n)\mid \Des_A(g)=\emptyset\}$.
One can show that $G(r,n)=US_n$, with $U\cap S_n=1$, i.e. for each $g\in G(r,n)$ there exist unique $\tau\in U$ and $\sigma \in S_n$ such that $g=\tau \sigma$ (see \cite[Proposition 4.1]{BZ1}). Moreover we have $\col(g)=\col(\tau)$ and $\Des_A(g)=\Des_A(\sigma)$. It follows that
\begin{eqnarray}
\nonumber\sum_{g\in G(r,n) }\sign(|g|)q^{\fmaj(g)}&=&\sum_{\tau\in U}\sum_{\sigma\in S_n}\sign(|\tau \sigma|)q^{\col(\tau)+r\cdot\maj(\sigma)}\\
\nonumber&=&\sum_{\tau\in U}\sign(|\tau|)q^{\col(\tau)}\, \sum_{\sigma\in S_n}\sign(\sigma)q^{r\cdot\maj(\sigma)}\\
\label{calc}&=&\sum_{\tau\in U}\sign(|\tau|)q^{\col(\tau)} [1]_{q^r}[2]_{-{q^r}}[3]_{q^r}\cdots [n]_{(-1)^{n-1}q^r},
\end{eqnarray}
by Theorem \ref{gessim}. Now we claim that 
\begin{equation}\label{U}
\sum_{\tau\in U}\sign(|\tau|)q^{\col(\tau)}=\left\{ \begin{array}{ll}[r]_{q^2}^m, & \textrm{if $n=2m$ is even;}\\  
{}[r]_q[r]_{q^2} ^{m},& \textrm{if $n=2m+1$ is odd.}\end{array}\right.
\end{equation}
In fact let $U(k_{r-1},\ldots,k_0)=\{\tau\in U \mid c_j(\tau)=i \textrm{ if }k_{r-1}+\cdots+k_{i+1}< j \leq k_{r-1}+\cdots+k_{i}\}$. In other words the elements in  $U(k_{r-1},\ldots,k_0)$ are those in $U$ having the first $k_{r-1}$ entries colored with $r-1$, then $k_{r-2}$ entries colored with $r-2$ and so on up to the last $k_0$ entries which are colored with $0$. 

It is clear that the map $\tau \mapsto |\tau|$ is a bijection between $U(k_{r-1},\ldots,k_0)$ and $S(k_{r-1},\ldots,k_0)$. Then
\begin{eqnarray*}
\sum_{\tau\in U}\sign(|\tau|)q^{\col(\tau)}&=& \sum_{k_{r-1},\ldots,k_0}\sum_{\tau\in U(k_{r-1},\ldots,k_0)} \sign(|\tau|)q^{\col(\tau)}\\
&=& \sum_{k_{r-1},\ldots,k_0} q^{\sum ik_i} \sum_{\sigma\in S(k_{r-1},\ldots,k_0)}\sign(\sigma)\\
&=& \sum_{k_{r-1},\ldots,k_0}\binom{\lfloor \frac{n}{2}\rfloor}{\lfloor\frac{k_{r-1}}{2}\rfloor ,\ldots,\lfloor \frac{k_{0}}{2}\rfloor}q^{\sum ik_i},
\end{eqnarray*}
by Lemma \ref{bin}, where the last sum is taken over all $r$-tuples $k_{r-1},\ldots,k_0$ such that $k_{r-1}+\cdots+k_0=n$ and at most one of the $k_i$'s is odd.
So, if $n=2m$ is even the $k_i$'s are all even, we let $k_i=2t_i$, and
\begin{eqnarray*}\sum_{\tau\in U}\sign(|\tau|)q^{\col(\tau)}&=&\sum_{\substack{t_{r-1},\ldots,t_0\\t_{r-1}+\cdots+t_0=m}}\binom{m}{t_{r-1},\ldots, t_0} q^{\sum 2it_i}\\
&=&(1+q^2+\cdots+q^{2(r-1)})^m\\
&=& [r]_{q^2}^m.  
\end{eqnarray*}
If $n=2m+1$ is odd there is exactly one of the $k_i$ which is odd. We split the sum according to this odd entry: if $k_i$ is odd we let $k_i=2t_i+1$ and if it is even we let $k_i=2t_i$. Then
\begin{eqnarray*}
\sum_{\tau\in U}\sign(|\tau|)q^{\col(\tau)}&=&\sum_{j=0}^{r-1}\sum_{\substack{t_{r-1},\ldots,t_0 \\t_{r-1}+\cdots+t_0=m}}\binom{m}{t_{r-1},\ldots ,t_0} q^{j+\sum 2it_i}\\
&=& \sum_{j=0}^{r-1}q^j \sum_{\substack{t_{r-1},\ldots,t_0\\t_{r-1}+\cdots+t_0=m}}\binom{m}{t_{r-1},\ldots, t_0} q^{\sum 2it_i}\\
&=& [r]_q[r]_{q^2}^m.
\end{eqnarray*}

Now one can easily check that 
\begin{equation}\label{rq2}
 [r]_{q^2}[2i-1]_{q^r}[2i]_{-q^r}=[(2i-1)r]_q[2ir]_{-q}.
\end{equation}
If $n$ is even the result follows immediately from Equations \eqref{calc}, \eqref{U}, and \eqref{rq2}. If $n$ is odd the result follows similarly using the further observation that $[r]_q[2m+1]_{q^r}=[(2m+1)r]_q$.
\end{proof}

Now we face the problem of the quotient groups $G(r,1,s,n)$. We need a preliminary result.
\begin{lem}\label{quoz}
Let $g\in G(r,1,s,n)$ and $g_0,\ldots,g_{s-1}$ be the $s$ lifts of $g$ in $G(r,n)$. We have

$$
\sum_{j=0}^{s-1}t^{\lambda_1(g_j)}q^{\fmaj(g_j)}=t^{\lambda_1(g)}q^{\fmaj(g)}[s]_{t^{r/s}q^{nr/s}}.$$ 

\end{lem}
\begin{proof}
We let $g_0$ be the unique lift of $g$ satisfying $c_n(g_{0})<r/s$. Then the other lifts are given by $|g_j|=|g_0|$ and $c_i(g_j)=R_r(c_i(g_0)+jr/s)$ for all $j\in[s-1]$ and for all $i\in[n]$. To compute $\fmaj$ we use definition \eqref{def-pezzi} with $s=1$, since all lifts $g_0,\ldots,g_{s-1} \in G(r,n)=G(r,1,1,n)$. It is clear that $\HDes(g_j)=\HDes(g_0)$ and $k_i(g_j)=k_i(g_0)+jr/s$. Therefore $\fmaj(g_j)=\fmaj(g_0)+njr/s$ and $\lambda_1(g_j)=\lambda_1(g)+jr/s$, so
\begin{eqnarray*}
\sum_{j=0}^{s-1}t^{\lambda_1(g_j)}q^{\fmaj(g_j)}&=&\sum_{j=0}^{s-1}t^{\lambda_1(g_0)+jr/s}q^{\fmaj(g_0)+jnr/s}\\
&=& t^{\lambda_1(g_0)}q^{\fmaj(g_0)}(1+t^{r/s}q^{nr/s}+\cdots+t^{(s-1)r/s}q^{(s-1)nr/s})=q^{\fmaj(g_0)}[s]_{t^{r/s}q^{nr/s}},
\end{eqnarray*}
and the result follows since $\fmaj(g_0)=\fmaj(g)$ and $\lambda_1(g_0)=\lambda_1(g)$.
\end{proof}
\begin{cor}\label{quotienttwist}We have
$$ \sum_{g\in G(r,1,s,n)}\sign(|g|)q^{\fmaj(g)}=[r]_q[2r]_{-q}[3r]_q\cdots [(n-1)r]_{(-1)^{n-2}q}\cdot [nr/s]_{(-1)^{n-1}q}.$$
\end{cor}
\begin{proof}
By Theorem \ref{grn} and Lemma \ref{quoz} specialized to the case $t=1$ we have
\begin{eqnarray*}
 \sum_{g\in G(r,1,s,n)}\sign(|g|)q^{\fmaj(g)}&=&\sum_{g\in G(r,n)}\sign(|g|)q^{\fmaj(g)}/[s]_{q^{nr/s}}\\
&=& [r]_q[2r]_{-q}[3r]_q\cdots [nr]_{(-1)^{n-1}q}/[s]_{q^{nr/s}},
\end{eqnarray*}
and the result easily follows since 
\begin{eqnarray*}
\frac{[nr]_{(-1)^{n-1}q}}{[s]_{q^{nr/s}}}&=&\frac{1-((-1)^{n-1}q)^{nr}}{1-(-1)^{n-1}q}\cdot \frac{1-q^{nr/s}}{1-q^{nr}}
= \frac{1-q^{nr}}{1-(-1)^{n-1}q}\cdot \frac{1-q^{nr/s}}{1-q^{nr}}\\
&=&\frac{1-((-1)^{n-1}q)^{nr/s}}{1-(-1)^{n-1}q}
=[nr/s]_{(-1)^{n-1}q},
\end{eqnarray*}
where we have used the fact that $r/s$ is an integer and hence $(n-1)nr/s$ is even.
\end{proof}

\begin{thm}\label{projtwist}
We have
$$
\sum_{g\in G(r,p,s,n)}\epsilon^{\inv(|g|)}q^{\fmaj(g)}=\left[\frac{r}{p}\right]_{ q^p}\left[\frac{2r}{p}\right]_{(\epsilon q)^p}\cdots \left[\frac{(n-1)r}{p}\right]_{(\epsilon^{n-2}q)^p}  \left[\frac{nr}{ps}\right]_{(\epsilon^{n-1}  q)^{p}}\Big\{[p]_{ q}^{n-m}[p]_{\epsilon  q}^{m} \Big\}_{q^p},
$$
where $m=\lfloor n/2 \rfloor$.
\end{thm}

\begin{proof}
It is a consequence of \cite{BB,Ca1} that 
$$\sum_{g\in G(r,1,s,n)}q^{\fmaj(g)}=[r]_q[2r]_{q}[3r]_q\cdots [(n-1)r]_{q}\cdot [nr/s]_{q},$$
and therefore we may unify this equation and Corollary~\ref{quotienttwist} by obtaining 
\begin{equation}\label{p0}
\sum_{g\in G(r,1,s,n)}\epsilon^{\inv(|g|)}q^{\fmaj(g)}=[r]_{ q}[2r]_{\epsilon q}[3r]_{ \epsilon^2q}\cdots [(n-1)r]_{{ \epsilon^{n-2}q}}\cdot [nr/s]_{{\epsilon^{n-1}q}},
\end{equation}
where $\epsilon=\pm 1$. We now want to describe the polynomial
$\sum_{g\in G(r,p,s,n)}\epsilon^{\inv(|g|)}q^{\fmaj(g)}$.
To this end we observe that an element $g\in G(r,1,s,n)$ belongs to $G(r,p,s,n)$ if and only if $\fmaj(g) \equiv 0 \mod p$ (see  \cite[Lemma 5.2]{Ca1}). Therefore we have
$$
\sum_{g\in G(r,p,s,n)}\epsilon^{\inv(|g|)}q^{\fmaj(g)}=\Bigg\{\sum_{g\in G(r,1,s,n)}\epsilon^{\inv(|g|)}q^{\fmaj(g)}\Bigg\}_{q^p}.
$$
We observe that
$$
\hspace{-7cm}[r]_{q}[2r]_{\epsilon q}[3r]_{\epsilon^{2}q}\cdots [(n-1)r]_{\epsilon^{n-2}q}\cdot [nr/s]_{\epsilon^{n-1}q}=
$$
\begin{eqnarray}\label{p1}
&=&\frac{1-q^r}{1-{q}} \cdot\frac{1-(\epsilon q)^{2r}}{1-{\epsilon q}}\cdots\frac{1-(\epsilon^{n-2}q)^{(n-1)r}}{1-{\epsilon^{n-2}q}}\cdots \frac{1-(\epsilon^{n-1}q)^{nr/s}}{1-{\epsilon^{n-1}q}} \nonumber \\
&=& \frac{(1-q^r)(1-q^{2r})\cdots (1-q^{(n-1)r})(1-q^{nr/s})}{(1-{q})(1-{\epsilon q})\cdots(1-{\epsilon^{n-1}q})}.
\end{eqnarray}

Notice that $\displaystyle{\big\{F(q)\big\}_{q^p}=\frac{1}{p}\sum_{j=0}^{p-1}F(\zeta^j_pq)}$. Hence from \eqref{p0} and \eqref{p1} we have
\begin{eqnarray*}
 \sum_{g\in G(r,p,s,n)}\epsilon^{\inv(|g|)}q^{\fmaj(g)}&=&\Big\{ [r]_{ q}[2r]_{ \epsilon q}[3r]_{ \epsilon^2q}\cdots [(n-1)r]_{{ \epsilon^{n-2}q}}\cdot [nr/s]_{{ \epsilon^{n-1}q}}\Big\}_{q^p}\\
&=& \frac{1}{p}\sum_{j=0}^{p-1}\frac{(1-(\zeta_p^jq)^r)(1-(\zeta_p^jq)^{2r})\cdots (1-(\zeta_p^jq)^{(n-1)r})(1-(\zeta_p^jq)^{nr/s})}{(1-{\zeta_p^jq})(1-{\epsilon \zeta_p^jq})\cdots(1-{\epsilon^{n-1}\zeta_p^jq})}\\
&=& \frac{1}{p}\sum_{j=0}^{p-1}\frac{(1-q^r)(1-q^{2r})\cdots (1-q^{(n-1)r})(1-q^{nr/s})}{(1-{\zeta_p^jq})(1-{\epsilon \zeta_p^jq})\cdots(1-{\epsilon^{n-1}\zeta_p^jq})}\\
&=&(1-q^r)(1-q^{2r})\cdots (1-q^{(n-1)r})(1-q^{nr/s})\\
&\cdot &  \frac{1}{p}\sum_{j=0}^{p-1}\frac{1}{(1-{\zeta_p^jq})(1-{\epsilon \zeta_p^jq})\cdots(1-{\epsilon^{n-1}\zeta_p^jq})}.
\end{eqnarray*}

The last factor can be manipulated as follows:

\begin{eqnarray*}
\frac{1}{p}\sum_{j=0}^{p-1}\frac{1}{(1-{\zeta_p^jq})(1-{\epsilon \zeta_p^jq})\cdots(1-{\epsilon^{n-1}\zeta_p^jq})}&=& \frac{1}{p}\sum_{j=0}^{p-1}\frac{[p]_{  \zeta_p^jq}}{(1- q^p)}\frac{[p]_{ \epsilon \zeta_p^jq}}{(1-(\epsilon q)^p)}\cdots \frac{[p]_{ \epsilon^{n-1} \zeta_p^jq}}{(1-(\epsilon^{n-1} q)^p)}\\
&=& \frac{\Big\{[p]_{q}^{n-m}[p]_{\epsilon q}^{m}\Big \}_{q^p}}{(1-q^p)(1-(\epsilon q)^p)\cdots(1-(\epsilon^{n-1} q)^p)},
\end{eqnarray*}

where $m=n/2$ if $n$ is even and $m=(n-1)/2$ if $n$ is odd. Hence we obtain that 

\begin{eqnarray*}
\sum_{g\in G(r,p,s,n)}\epsilon^{\inv(|g|)}q^{\fmaj(g)} &=& \frac{(1-q^r)(1-q^{2r})\cdots (1-q^{(n-1)r})(1-q^{nr/s})}{(1- q^p)(1-(\epsilon q)^p)\cdots(1-(\epsilon^{n-1} q)^p)}  \Big\{[p]_{q}^{n-m}[p]_{\epsilon q}^{m}\Big \}_{q^p}\\
&=& \frac{1-q^r}{1-q^p}\cdots \frac{1-q^{(n-1)r}}{1-(\epsilon^{n-2} q)^p}\frac{1-q^{nr/s}}{1-(\epsilon^{n-1} q)^p}\Big\{[p]_{q}^{n-m}[p]_{\epsilon q}^{m}\Big \}_{q^p}\\
&=& \left[\frac{r}{p}\right]_{q^p}\left[\frac{2r}{p}\right]_{(\epsilon q)^p}\cdots \left[\frac{(n-1)r}{p}\right]_{(\epsilon^{n-2}q)^p}  \left[\frac{nr}{ps}\right]_{(\epsilon^{n-1}  q)^{p}}\Big\{[p]_{q}^{n-m}[p]_{\epsilon q}^{m} \Big\}_{q^p}.
\end{eqnarray*}
\end{proof}

Now we are ready to prove our main result.

\begin{proof}[{\it Proof of Theorem~\ref{maintwist}.}]
Now Theorem \ref{maintwist} is  a direct consequence of Theorem \ref{projtwist} since $k\cdot\fmaj(g)\equiv k\cdot c(g)\mod r$ for all $g\in \qc$ and so 
$$
\sum_{g\in G(r,p,s,n)}\epsilon^{\inv(|g|)}\zeta_r^{k\cdot c(g)}q^{\fmaj(g)}=\sum_{g\in G(r,p,s,n)}\epsilon^{\inv(|g|)}(\zeta_r^k q)^{\fmaj(g)}
$$
\end{proof}
\section{Carlitz's Identities}
In this section we give a general method to compute the trivariate distribution of $\des$ (or $\fdes$), $\fmaj$ and $\col$ over $\qc$. This unifies and generalizes all related results cited in the introduction, and provides two different generalizations of Carlitz's identity for the group $\qc$.
\smallskip

For $f=(f_1,\ldots,f_n)\in \mathbb N^n$ let $|f|=f_1+\cdots+f_n$ and $$
\mathbb N^n(p)=\{f\in \mathbb N^n \mid |f| \equiv 0 \mod p\}.
$$
Moreover, if $f\in \mathbb N^n$, for any $k >0$  we let 
$$\col_k(f)\eqdef\sum_{i=1}^nR_{k}(f_i)$$ be the sum of the residues of the $f_i$'s modulo $k$.

\begin{lem}\label{bijcar}
There is a bijection
$$
\mathbb N^n(p)\longleftrightarrow G(r,p,s,n) \times \mathcal P_n \times [0,s-1],
$$
such that if $f\leftrightarrow (g,\lambda,h)$ then
\begin{enumerate}
\item[$(a)$] $\max f=\lambda_1(g)+r\lambda_1+hr/s$;
\item[$(b)$] $|f|=\fmaj(g)+r|\lambda|+hnr/s$;
\item[$(c)$] $\col_{\small{\frac{r}{s}}}(f)=\col(g)$.
\end{enumerate}
\end{lem}
\begin{proof}
This result is similar to several others appearing in the literature and in particular is a special case of the bijection appearing in \cite[Theorem 8.3]{Ca1} and so we simply describe how the bijection is defined for the reader's convenience. 
If $f=(f_1,\ldots,f_n)\in \mathbb N^n(p)$, then $g$ is the unique element in $\qc$ having a lift $\tilde g\in G(r,p,n)$ satisfying:
\begin{itemize}
\item $f_{|g(i)|}\geq f_{|g(i+1)|}$ for all $i\in [n-1]$;
\item if $f_{|g(i)|}=f_{|g(i+1)|}$ then $|g(i)|<|g(i+1)|$;
\item $c_i(\tilde g)\equiv f_{|g(i)|}\mod r$ for all $i\in [n]$.
\end{itemize}
Letting $\mu$ be the partition obtained by reordering the entries in $(f_1,\ldots,f_n)$ in nonincreasing order one can show that 
$\mu-\lambda(g)$ (componentwise difference) is still a partition whose parts are all congruent to the same multiple of $r/s \mod r$ (i.e. $\mu$ is $g$-compatible in the notation of \cite{Ca1}). The partition $\lambda\in \mathcal P_n$ and the integer $h\in [0,s-1]$ are therefore uniquely determined by the requirement that $\mu=\lambda(g)+r\cdot \lambda+(hr/s,\ldots,hr/s)$.

The inverse map of this bijection is much simpler. Let $(g,\lambda,h)$ be a triple in $\qc\times\mathcal P_n \times [0,s-1]$, then the corresponding element $f\in \mathbb N^n(p)$ is given by $f_i=\lambda_{|g^{-1}(i)|}(g)+r\lambda_{|g^{-1}|(i)}+hr/s$.

All the other statements are straightforward consequences, since $\lambda_i(g)\equiv c_i(\tilde g)\mod r/s$ for any lift $\tilde g$ of $g$, and so also $\col(g)=\sum_{i=1}^n R_{r/s}(\lambda_i(g))$.
\end{proof}
We will make use of the following observation several times: if $\mathcal F$ is a set and $m:\mathcal F\rightarrow \mathbb N$ is such that the set $\{f\in \mathcal F \mid m(f)\leq k\}$ is finite for all $k\in \mathbb N$ then we always have
\begin{equation}\label{leqeq}
 \frac{1}{1-t}\sum_{f\in \mathcal F}t^{m(f)}\phi(f)=\frac{1}{1-t}\sum_{k\in \mathbb N}t^k\sum_{f\in \mathcal F \mid m(f)= k}\phi(f)=\sum_{k\in \mathbb N}t^k\sum_{f\in \mathcal F \mid m(f)\leq k}\phi(f),
\end{equation}
where $\phi$ is any map defined on $\mathcal F$ taking values in a (polynomial) ring.

We can now state the main result of this section.
\begin{thm}\label{desfmajcol}
$$
\Bigg\{\sum_{k\geq 0}t^k \left ([k+1]_{q^{r/s}}+aq[k]_{q^{r/s}}\Big[\frac{r}{s}-1\Big]_{aq}\right)^n\Bigg\}_{q^p}=\frac{\sum_{g\in G(r,p,s,n)}t^{\des(g)}q^{\fmaj(g)}a^{\col(g)}}{(1-t)(1-t^sq^r)\cdots(1-t^sq^{(n-1)r})(1-tq^{nr/s})}
$$
\end{thm}
\begin{proof}
 For $f\in \mathbb N^n(p)$ we let $m(f)=\lfloor \frac{s \max (f)+r-s}{r}\rfloor$ and we observe that $m(f)\leq k$ if and only if $\max(f)\leq \frac{r}{s}k$. In fact $\lfloor\frac{s \max (f)+r-s}{r}\rfloor\leq k$ is equivalent to $s \max (f)+r-s\leq kr+r-1$. This happens if and only if $\max (f)\leq kr/s+\frac{s-1}{s}$ which is equivalent to $\max (f)\leq kr/s$ since $\max(f)$ is an integer.
 We want to make use of Equation \eqref{leqeq} with $\phi(f)=q^{|f|}a^{\col_{\small{\frac{r}{s}}}(f)}$ and so we compute

\begin{eqnarray*}
\sum_{\substack{f\in \mathbb N^n(p)\mid \\ m(f)\leq k}}q^{|f|}a^{\col_{\small{\frac{r}{s}}}(f)}&=
&\bigg\{\Big(\sum_{j=0}^{kr/s}q^ja^{R_{r/s}(j)}\Big)^n\bigg\}_{q^p}\\
&=& \bigg\{\left(1+q^{\frac{r}{s}}+\cdots+q^{\frac{kr}{s}}+(aq+a^2q^2+\cdots+a^{\frac{r}{s}-1}q^{\frac{r}{s}-1})(1+q^{\frac{r}{s}}+\cdots+q^{(k-1)\frac{r}{s}})\right)^n\bigg\}_{q^p}
\\&=&
\bigg\{\Big ([k+1]_{q^{r/s}}+aq\Big[\frac{r}{s}-1\Big]_{aq}[k]_{q^{r/s}}\Big)^n\bigg\}_{q^p}.
\end{eqnarray*}
Therefore, by Equation \eqref{leqeq},
\begin{eqnarray*}
\sum_{f\in \mathbb N^n(p)}t^{m(f)}q^{|f|}a^{\col_{\small{\frac{r}{s}}}(f)}&=&(1-t)\Bigg\{\sum_{k\geq 0}t^k \bigg ([k+1]_{q^{r/s}}+aq\Big[\frac{r}{s}-1\Big]_{aq}[k]_{q^{r/s}}\bigg)^n\Bigg\}_{q^p}.
\end{eqnarray*}
Now we want to compute this polynomial using the bijection described in Lemma \ref{bijcar}. Observe that if $f\mapsto(g,\lambda,h)$ then $m(f)=\des(g)+s\lambda_1+h$. Therefore we have
\begin{eqnarray*}
\sum_{f\in \mathbb N^n(p)}t^{m(f)}q^{|f|}a^{\col_{\small{\frac{r}{s}}}(f)}&=&\sum_{(g,\lambda,h)}t^{\des(g)+s\lambda_1+h}q^{\fmaj(g)+r|\lambda|+nhr/s}a^{\col(g)}\\
&=& \sum_{g\in G(r,p,s,n)}t^{\des(g)}q^{\fmaj(g)}a^{\col(g)}\sum_{\lambda\in \mathcal P_n}t^{s\lambda_1}q^{r|\lambda|}\sum_{h=0}^{s-1}t^{h}q^{nhr/s}\\
&=& \frac{\sum_{g\in G(r,p,s,n)}t^{\des(g)}q^{\fmaj(g)}a^{\col(g)}}{(1-t^sq^r)\cdots(1-t^sq^{(n-1)r})(1-tq^{nr/s})},
\end{eqnarray*} 
and the result follows.
\end{proof}
Letting $s=p=1$ in the previous result we obtain \cite[Equation (8.1)]{BZ1}. Moreover,  for $a=1$ we have $[k+1]_{q^{r/s}}+q[k]_{q^{r/s}}[\frac{r}{s}-1]_q=\Big[\frac{r}{s}k+1\Big]_q$ and hence we obtain the following result.
\begin{cor}[Carlitz's identity for $\qc$ with $\des$]\label{cardes}
$$\Bigg\{\sum_{k\geq0} t^k \Big[\frac{r}{s}k+1\Big]_q^n\Bigg\}_p=\frac{\sum_{g\in G(r,p,s,n)}t^{\des(g)}q^{\fmaj(g)}}{(1-t)(1-t^sq^r)\cdots(1-t^sq^{(n-1)r})(1-tq^{nr/s})}.
$$\end{cor}
The special case with $r=2$, $p=s=1$ of Corollary \ref{cardes} is the main result of \cite{CG}, and for $p=s=1$ we obtain \cite[Theorem 10 (iv)]{CM}.

Now we see how a simple modification of the same ideas lead to the generalization of other identities that use different flag-descents.
\begin{thm}
\begin{eqnarray}\label{long}
\sum_{k\geq 0} t^k \left( [Q_{r/s}(k)+1]_{q^{r/s}} +aq[r/s-1]_{aq} \cdot [Q_{r/s}(k)]_{q^{r/s}} + aq^{mr/s+1)}[R_{r/s}(k)]_{aq}\right)^n\\
=\frac{\sum_{g\in G(r,p,s,n)} t^{\fdes(g)}q^{\fmaj(g)} a^{\col(g)}} {(1-t)(1-t^s q^r)\cdots (1-t^s q^{(n-1)r/s})\cdots (1-t q^{nr/s})},
\end{eqnarray}
where $Q_{r/s}(k)$ is determined by $k=r/s \cdot Q_{r/s}(k)+ R_{r/s}(k)$.
\end{thm}

A proof of this result can be done by paralleling that of  Theorem \ref{desfmajcol}. In this case one should simply consider the statistics $m(f)=\max(f)$ and follow the rest mutatis mutandis.
By letting $a=1$ in Equation~\ref{long}, one obtain a second Carlitz's identity type, with flag-descents.

%
%
%

\begin{thm}[Carlitz's identity of $G(r,p,s,n)$ with $\fdes$]\label{Carlitz-ps}
$$
\bigg\{ \sum_{k\geq 0}t^k [k+1]^n_q\bigg\}_{q^p}=\frac{\sum_{g\in G(r,p,s,n)} t^{\fdes(g)}q^{\fmaj(g)}}{(1-t)(1-t^r q^r)(1-t^r q^{2r})\cdots (1-t^r q^{(n-1)r}) (1-t^\frac{r}{s} q^{\frac{nr}{s}})}.
$$
\end{thm}
For $r=2$ and $p=s=1$ we obtain \cite[Theorem 4.2]{ABR1}, for $r=s=2$ and $p=1$ \cite[Theorem 4.3]{BC}, and for $p=1$ \cite[Theorem 11.2]{BB}.


\section{Multivariate generating functions}\label{genefunc}
In this section we make further use of the bijection \cite[Theorem 8.3]{Ca1} to compute new multivariate distributions for the groups $\qc$. We first concentrate on the case of the groups $G(r,n)$ to make the used arguments more clear to the reader. We need to state the particular case of this bijection needed for our purposes, namely the special case of 2-partite partitions.
\smallskip

We recall that a \emph{$2$-partite partition of length $n$} (see \cite{GG}) is a $2\times n$ matrix with non negative integer coefficients  $f=\left(\begin{array}{cccc}f^{(1)}_1&f^{(1)}_2&\ldots&f^{(1)}_n\\f^{(2)}_1&f^{(2)}_2&\ldots &f^{(2)}_n\end{array}\right)$ satisfying the following conditions:
\begin{itemize}
\item $f^{(1)}_1\geq f^{(1)}_2\geq\ldots \geq f^{(1)}_n$;
\item If $f^{(1)}_i = f^{(1)}_{i+1}$ then $f^{(2)}_i \geq f^{(2)}_{i+1}$.
\end{itemize}

One may think of a 2-partite partition as a generic multiset of pairs of nonnegative integers (the columns of $f$) of cardinality $n$.
If $f$ is a 2-partite partition we denote by $f^{(1)}$ and $f^{(2)}$ the first and the second row of $f$ respectively.
To be consistent with the notation of \cite[\S 8]{Ca1}, we denote by $\mathcal B(n)$ the set of $2$-partite partitions of length $n$ and we let
\begin{align*}
\mathcal B(r,n)&\eqdef\{f\in \mathcal B(n) \mid  f^{(1)}_i+f^{(2)}_i\equiv 0 \mod r \textrm { for all }i\in[n]\} \qquad {\rm and} \\
\mathcal B(r,s,1,n)&\eqdef\{f\in \mathcal B(n) \mid \textrm{ there exists } l\in [s-1] \textrm{ such that }f^{(1)}_i+f^{(2)}_i\equiv lr/s \mod r \textrm { for all }i\in[n]\}.
\end{align*}
Note that $\mathcal B(1,n)=\mathcal B(n)$ and that $\mathcal B(r,n)=\mathcal B(r,1,1,n)$.

\begin{prop}\label{bijGrn}There exists a bijection between $\mathcal B(r,s,1,n)$ and the set of the 5-tuples $(g,\lambda,\mu,h,k)$, where $g\in G(r,1,s,n)$, $\lambda, \mu\in \mathcal P_n$, and $h,k\in [0,s-1]$.
In this bijection, if $f\leftrightarrow(g,\lambda,\mu,h,k)$ then
\begin{enumerate}
\item[$(a)$] $\max(f^{(1)})=\lambda_1(g)+r\lambda_1+hr/s$;
\item[$(b)$] $\max(f^{(2)})=\lambda_1(g^{-1})+r\mu_1+kr/s$;
\item[$(c)$] $|f^{(1)}|=\fmaj(g)+r|\lambda|+hnr/s$;
\item[$(d)$] $|f^{(2)}|=\fmaj(g^{-1})+r|\mu|+knr/s$;
\item[$(e)$] $\col_{\small{\frac{r}{s}}}(f^{(1)})= \col(g) $;
\item[$(f)$] $\col_{\small{\frac{r}{s}}}(f^{(2)})=\col(g^{-1})$.
\end{enumerate} 
\end{prop}

\begin{proof}
This is again a particular case of \cite[Theorem 8.3]{Ca1}. In this case the bijection is defined as follows: if $(g,\lambda,\mu,h,k)\leftrightarrow f$ then
$$
f^{(1)}_i=\lambda_i(g)+r\lambda_i+hr/s\textrm{ and } f^{(2)}_i=\lambda_{|g(i)|}(g^{-1})+r\mu_{|g(i)|}+kr/s,
$$
and all the other statements follow immediately.
\end{proof}

For the reader's convenience we first face the problem of the wreath products $G(r,n)$.
\begin{thm}\label{sixgrn} Let $r \in \NN$. Then 
$$
\hspace{-6cm}\sum_{k_1,k_2\geq 0}t_1^{k_1}t_2^{k_2}\bigg(\prod_{\substack{i\in[0,rk_1],\, j\in [0,rk_2]:\\i+j\equiv 0 \!\mod r}}\frac{1}{1-ua_1^{R_r(i)}a_2^{R_r(j)}q_1^{i}q_2^{j}}\bigg)=$$
$$=\sum_{n\geq 0}u^n\sum_{g\in G(r,n)}t_1^{\des(g)}t_2^{\des(g^{-1})}q_1^{\fmaj(g)}, q_2^{\fmaj(g^{-1})}a_1^{\col(g)}a_2^{\col(g^{-1})}\prod_{j=0}^{n}\frac{1}{(1-t_1q_1^{jr})(1-t_2q_2^{jr})}.
$$
\end{thm}

\begin{proof} For $i=1,2$, we let $m_i(f)=\lfloor \frac{\max f^{(i)}+r-1}{r}\rfloor$
and $\mathcal B_{k_1,k_2}=\{f\in \mathcal B(r,n) \mid m_i(f)\leq k_i, \ i=1,2\}$. We can easily verify that $m_i(f)\leq k_i$ if and only if $\max(f^{(i)})\leq rk_i$. Therefore, using an idea already appearing in \cite{GG}, one can show that
\begin{eqnarray*}
\sum_{f\in \mathcal B_{k_1,k_2}}q_1^{|f^{(1)}|} q_2^{|f^{(2)}|}a_1^{\col_{\small{\frac{r}{s}}}(f^{(1)})}a_2^{\col_{\small{\frac{r}{s}}}(f^{(2)})}
&=& \bigg(\prod_{\substack{i\in[0,rk_1],\, j\in [0,rk_2]:\\i+j\equiv 0 \!\mod r}}\frac{1}{1-ua_1^{R_r(i)}a_2^{R_r(j)}q_1^{i}q_2^{j}}\bigg)\bigg|_{u^n}.
\end{eqnarray*}
\newpage

Now we use a natural bivariate extension of Equation \eqref{leqeq} to obtain
\begin{eqnarray*}
\hspace{-7cm} \sum_{k_1,k_2 \geq 0}t_1^{k_1}t_2^{k_2}\sum_{f\in \mathcal B_{k_1,k_2}}q_1^{|f^{(1)}|} q_2^{|f^{(2)}|}a_1^{\col_{\small{\frac{r}{s}}}(f^{(1)})}a_2^{\col_{\small{\frac{r}{s}}}(f^{(2)})}
\end{eqnarray*}
\begin{eqnarray*}
\hspace{3cm}=\frac{1}{(1-t_1)(1-t_2)}\sum_{f\in \mathcal B(r,n)}t_1^{m_1(f)}t_2^{m_2(f)}q_1^{|f^{(1)}|} q_2^{|f^{(2)}|}a_1^{\col_{\small{\frac{r}{s}}}(f^{(1)})}a_2^{\col_{\small{\frac{r}{s}}}(f^{(2)})}
\end{eqnarray*}
and hence
$$
\hspace{-7cm} \sum_{f\in \mathcal B(r,n)}t_1^{m_1(f)}t_2^{m_2(f)}q_1^{|f^{(1)}|} q_2^{|f^{(2)}|}a_1^{\col_{\small{\frac{r}{s}}}(f^{(1)})}a_2^{\col_{\small{\frac{r}{s}}}(f^{(2)})}=
$$
\begin{eqnarray*}&=&(1-t_1)(1-t_2) \sum_{k_1,k_2\geq 0}t_1^{k_1}t_2^{k_2}\sum_{f\in \mathcal B_{k_1,k_2}}q_1^{|f^{(1)}|} q_2^{|f^{(2)}|}a_1^{\col_{\small{\frac{r}{s}}}(f^{(1)})}a_2^{\col_{\small{\frac{r}{s}}}(f^{(2)})}\\
&=& (1-t_1)(1-t_2)\sum_{k_1,k_2\geq 0 }t_1^{k_1}t_2^{k_2}\bigg(\prod_{\substack{i\in[0,rk_1],\, j\in [0,rk_2]:\\i+j\equiv 0 \!\mod r}}\frac{1}{1-ua_1^{R_r(i)}a_2^{R_r(j)}q_1^{i}q_2^{j}}\bigg)\bigg|_{u^n}.
\end{eqnarray*}
On the other hand, by Proposition \ref{bijGrn} used with $s=1$, we have

$$
\hspace{-7cm} \sum_{f\in \mathcal B(r,n)}t_1^{m_1(f)}t_2^{m_2(f)}q_1^{|f^{(1)}|} q_2^{|f^{(2)}|}a_1^{\col_{\small{\frac{r}{s}}}(f^{(1)})}a_2^{\col_{\small{\frac{r}{s}}}(f^{(2)})}=
$$
\begin{eqnarray*}
 &=&\sum_{(g,\lambda,\mu)}t_1^{\des(g)+\lambda_1}t_2^{\des(g^{-1})+\mu_1}q_1^{\fmaj(g)+r|\lambda|} q_2^{\fmaj(g^{-1})+r|\mu|}a_1^{\col(g)}a_2^{\col(g^{-1})}\\
&=& \sum_{g\in G(r,n)}t_1^{\des(g)}t_2^{\des(g^{-1})}q_1^{\fmaj(g)} q_2^{\fmaj(g^{-1})}a_1^{\col(g)}a_2^{\col(g^{-1})}
\sum_{\lambda \in \mathcal{P}_n}t_1^{\lambda_1}q_1^{r|\lambda|}\sum_{\mu \in \mathcal{P}_n}t_2^{\mu_1}q_2^{r|\mu|}\\
&=& \sum_{g\in G(r,n)}t_1^{\des(g)}t_2^{\des(g^{-1})}q_1^{\fmaj(g)} q_2^{\fmaj(g^{-1})}a_1^{\col(g)}a_2^{\col(g^{-1})}\prod_{j=1}^{n}\frac{1}{(1-t_1q_1^{jr})(1-t_2q_2^{jr})}
\end{eqnarray*}
and the proof is complete.
\end {proof}

Now we want to extend this result to the quotients $G(r,1,s,n)$. 
\begin{thm}\label{apap} Let $r,s \in \NN$, such that $s|r$. Then 
$$
\hspace{-6cm}\sum_{k_1,k_2\geq 0 }t_1^{k_1}t_2^{k_2}\sum_{l=0}^{s-1}\bigg(\prod_{\substack{i\in [0,k_1\frac{r}{s}],\,j\in [0,k_2\frac{r}{s}]:\\i+j\equiv l\frac{r}{s}\!\mod r}}\frac{1}{1-ua_1^{R_{r/s}(i)}a_2^{R_{r/s}(j)}q_1^{i}q_2^{j}}\bigg)=$$
$$\hspace{+2cm}=\sum_{n\geq 0} u^n \frac{\sum_{g\in G(r,1,s,n)}t_1^{\des(g)}t_2^{\des(g^{-1})}q_1^{\fmaj(g)} q_2^{\fmaj(g^{-1})}a_1^{\col(g)}a_2^{\col(g^{-1})}}{{(1-t_1)(1-t_2)(1-t_1q_1^{n\frac{r}{s}})(1-t_2q_2^{n\frac{r}{s}})}\prod_{j=1}^{n-1}{(1-t_1^sq_1^{jr})(1-t_2^sq_2^{jr})}}.
$$
\end{thm}

\begin{proof} 
Say that $f\in \mathcal B(r,s,1,n)$ is of type $l$, with $l\in [0,s-1]$, if the column sums of $f$ are all congruent to $lr/s \mod r$.
For $i=1,2$, we let $m_i(f)=\lfloor \frac{s \max (f^{(i)})+r-s}{r}\rfloor$ and $\mathcal B_{k_1,k_2}^{(l)}=\{f\in \mathcal B(r,s,1,n)\textrm{ of type $l$} \mid m_1(f)\leq k_1, m_2(f)\leq k_2) \}$. We observe that $m_i(f)\leq k_i$ if and only if $\max(f^{(i)})\leq \frac{r}{s}k_i$ and we compute
$$ \sum_{f\in \mathcal B_{k_1,k_2}^{(l)}}q_1^{|f^{(1)}|} q_2^{|f^{(2)}|}a_1^{\col_{\small{\frac{r}{s}}}(f^{(1)})}a_2^{\col_{\small{\frac{r}{s}}}(f^{(2)})}= \bigg(\prod_{\substack{i\in [0,k_1\frac{r}{s}],\,j\in [0,k_2\frac{r}{s}]:\\i+j\equiv l\frac{r}{s}\!\mod r}}\frac{1}{1-ua_1^{R_{r/s}(i)}a_2^{R_{r/s}(j)}q_1^{i}q_2^{j}}\bigg)\bigg|_{u^n}.$$

Therefore
$$\hspace{-5cm}\frac{1}{(1-t_1)(1-t_2)}\sum_{f\in \mathcal B(r,s,1,n)}t_1^{m_1(f)}t_2^{m_2(f)}q_1^{|f^{(1)}|} q_2^{|f^{(2)}|}a_1^{\col_{\small{\frac{r}{s}}}(f^{(1)})}a_2^{\col_{\small{\frac{r}{s}}}(f^{(2)})}$$
\begin{eqnarray*}
&=&\sum_{k_1,k_2\geq 0}t_1^{k_1}t_2^{k_2}\sum_{l\in [0,s-1]}\sum_{f\in \mathcal B^{(l)}_{k_1,k_2}}q_1^{|f^{(1)}|} q_2^{|f^{(2)}|}a_1^{\col_{\small{\frac{r}{s}}}(f^{(1)})}a_2^{\col_{\small{\frac{r}{s}}}(f^{(2)})}\\
&=&\sum_{k_1,k_2\geq 0}t_1^{k_1}t_2^{k_2}\sum_{l\in [0,s-1]}\bigg(\prod_{\substack{i\in [0,k_1\frac{r}{s}],\,j\in [0,k_2\frac{r}{s}]:\\i+j\equiv l\frac{r}{s}\!\mod r}}\frac{1}{1-ua_1^{R_{r/s}(i)}a_2^{R_{r/s}(j)}q_1^{i}q_2^{j}}\bigg)\bigg|_{u^n}.
\end{eqnarray*}

On the other hand, by Proposition \ref{bijGrn}, we have
$$
\hspace{-9cm} \sum_{f\in \mathcal B(r,s,1,n)}t_1^{m_1(f)}t_2^{m_2(f)}q_1^{|f^{(1)}|} q_2^{|f^{(2)}|}a_1^{\col_{\small{\frac{r}{s}}}(f^{(1)})}a_2^{\col_{\small{\frac{r}{s}}}(f^{(2)})}
$$
\begin{eqnarray*}
&=&\sum_{g,\lambda,\mu,h,k}t_1^{\des(g)+s\lambda_1+h}t_2^{\des(g^{-1})+s\mu_1+k}q_1^{\fmaj(g)+r|\lambda|+kn\frac{r}{s}} q_2^{\fmaj(g^{-1})+r|\mu|+hn\frac{r}{s}}a_1^{\col(g)}a_2^{\col(g^{-1})}\\
&=&\sum_{g\in G(r,1,s,n) }t_1^{\des(g)}t_2^{\des(g^{-1})}q_1^{\fmaj(g)}q_2^{\fmaj(g^{-1})}a_1^{\col(g)}a_2^{\col(g^{-1})}\sum_{h=0}^{s-1}t_1^h q_1^{hn\frac{r}{s}}\sum_{k=0}^{s-1}t_2^k q_2^{kn\frac{r}{s}}\sum_{\lambda } t_1^{s\lambda_1}q_1^{r|\lambda|}\sum_{\mu}t_2^{s\mu_1}q_2^{r|\mu|}\\
&=& \sum_{g\in G(r,1,s,n) }t_1^{\des(g)}t_2^{\des(g^{-1})}q_1^{\fmaj(g)}q_2^{\fmaj(g^{-1})}a_1^{\col(g)}a_2^{\col(g^{-1})}\frac{1-t_1^sq_1^{nr}}{1-t_1q_1^{n\frac{r}{s}}}\frac{1-t_2^sq_2^{nr}}{1-t_2q_2^{n\frac{r}{s}}}\prod_{j=1}^n \frac{1}{(1-t_1^sq_1^{jr})(1-t_2^sq_2^{jr})}.
\end{eqnarray*}
\end{proof}

The corresponding result for the the full class of projective reflection groups $G(r,p,s,n)$ can now be easily deduced from  Theorem \ref{apap}. We remind that $\qc$ is defined if and only if $p|r$, $s|r$, and $sp|rn$. In particular if $r,p,s$ are fixed then $\qc$ is defined only if $sp/ {\rm GCD}(sp,r)$ divides $n$. 

\begin{thm}\label{adg} Let $r,p,s \in \NN$, such that $s$ and $p$ divide $r$. Let $d= sp/ {\rm GCD}(sp,r)$. Then  
$$
\hspace{-5cm}\Bigg\{\sum_{k_1,k_2\geq 0}t_1^{k_1}t_2^{k_2}\sum_{l=0}^{s-1}\bigg(\prod_{\substack{i\in [0,k_1\frac{r}{s}],\,j\in [0,k_2\frac{r}{s}]:\\i+j\equiv l\frac{r}{s}\!\mod r}}\frac{1}{1-ua_1^{R_{r/s}(i)}a_2^{R_{r/s}(j)}q_1^{i}q_2^{j}}\bigg)\Bigg\}_{u^d q_1^p}=$$
$$=\sum_{n \geq 0 \; : \; d|n}u^n \frac{\sum_{g\in \qc}t_1^{\des(g)}t_2^{\des(g^{-1})}q_1^{\fmaj(g)} q_2^{\fmaj(g^{-1})}a_1^{\col(g)}a_2^{\col(g^{-1})}}{{(1-t_1)(1-t_2)(1-t_1q_1^{n\frac{r}{s}})(1-t_2q_2^{n\frac{r}{s}})}\prod_{j=1}^{n-1}{(1-t_1^sq_1^{jr})(1-t_2^sq_2^{jr})}}.
$$
\end{thm}

\subsection{Hilbert series of diagonal invariant algebras}

We can now make a specialization of this result. We set $a_1=a_2=1$ in Theorem~\ref{adg}
$$
\hspace{-5cm}\Bigg\{\sum_{k_1,k_2\geq 0}t_1^{k_1}t_2^{k_2}\sum_{l=0}^{s-1}\bigg(\prod_{\substack{i\in [0,k_1\frac{r}{s}],\,j\in [0,k_2\frac{r}{s}]:\\i+j\equiv l\frac{r}{s}\!\mod r}}\frac{1}{1-uq_1^{i}q_2^{j}}\bigg)\Bigg\}_{u^d q_1^p}=$$
$$=\sum_{n\geq 0  \; : \; d|n }u^n \frac{\sum_{g\in \qc}t_1^{\des(g)}t_2^{\des(g^{-1})}q_1^{\fmaj(g)} q_2^{\fmaj(g^{-1})}}{{(1-t_1)(1-t_2)(1-t_1q_1^{n\frac{r}{s}})(1-t_2q_2^{n\frac{r}{s}})}\prod_{j=1}^{n-1}{(1-q_1^{jr}t_1^s)(1-q_2^{jr}t_2^s)}}.
$$
Multiplying both sides of the previous identity by $(1-t_1)(1-t_2)$ and then taking the limit as $t_1,t_2 \rightarrow 1$ we obtain
\begin{equation}\label{e:hilb}
\Bigg\{\sum_{l=0}^{s-1}\prod_{\substack{i+j\equiv lr/s \\ \mod r}}\frac{1}{1-uq_1^{i}q_2^{j}}\Bigg\}_{u^d q_1^p}=\sum_{n\geq 0  \; : \; d|n}u^n\frac{\sum_{g\in \qc}q_1^{\fmaj(g)} q_2^{\fmaj(g^{-1})}} {(1-q_1^{n\frac{r}{s}})(1-q_2^{n\frac{r}{s}})\prod_{j=1}^{n-1}{(1-q_1^{jr})(1-q_2^{jr})}}.
\end{equation}
There is a nice algebraic interpretation of the previous identity that goes back to the work of Adin and Roichman \cite{AR}.  We let $S_p[X,Y]$ be the subalgebra of the algebra of polynomials in $2n$ variables $x_1,\ldots,x_n,y_1,\ldots,y_n$ generated by 1 and the monomials whose degrees in both the $x$'s and the $y$'s variables is divided by $p$. Then we can observe (see \cite[\S 8]{Ca1}) that the factor
$$
\frac{1}{(1-q_1^{n\frac{r}{s}})(1-q_2^{n\frac{r}{s}})\prod_{j=1}^{n-1}{(1-q_1^{jr})(1-q_2^{jr})}}
$$
is the bivariate Hilbert series of the invariant algebra of the ``tensorial action'' of the group $G(r,s,p,n)^2=G(r,s,p,n)\times G(r,s,p,n)$ (note the interchanging of the roles of $p$ and $s$) on the ring of polynomials $S_p[X,Y]$. Furthermore by \cite[Corollary 8.6]{Ca1} we can deduce that
\begin{equation}\label{hilbert}
\frac{\sum_{g\in \qc}q_1^{\fmaj(g)} q_2^{\fmaj(g^{-1})}} {(1-q_1^{n\frac{r}{s}})(1-q_2^{n\frac{r}{s}})\prod_{j=1}^{n-1}{(1-q_1^{jr})(1-q_2^{jr})}}=\Hilb(S_p[X,Y]^{\Delta G(r,s,p,n)}),
\end{equation}
where $\Delta G(r,s,p,n)$ denotes the diagonal embedding of $G(r,s,p,n)$ in $G(r,s,p,n)^2$. 
We can now conclude that Equation \eqref{e:hilb} provides the following interpretation for the generating function of the Hilbert series of the diagonal invariant algebras of the groups $\qc$.
\begin{cor}
Let $r,p,s \in \NN$, such that $s$ and $p$ divide $r$. Let $d= sp/ {\rm GCD}(sp,r)$. Then
$$
\sum_{n\geq 0 \; : \; d|n }u^n\Hilb \big(S_p[X,Y]^{\Delta \qc}\big)(q_1,q_2)=\Bigg\{\sum_{l=0}^{p-1}\prod_{\substack{i+j\equiv lr/s \\ \mod r}}\frac{1}{1-uq_1^{i}q_2^{j}}\Bigg\}_{u^dq_1^s}.
$$ 
\end{cor}
This identity was apparently known only in the case of the symmetric group \cite{GG,Ro}.

\section{Orderings}\label{orde}

In the literature, when considering statistics based on descents and/or inversions on the groups $G(r,n)$ at least two distinct orderings have been considered on the set of colored integers. In this last section we would like to clarify some of the relationships between two of these orderings concerning the results shown in this paper. 

We recall that we have been using the order defined in \eqref{color-order}.
One can consider also the order $<'$
$$
n^{r-1} <'\cdots<' n^1<' \cdots <' 1^{r-1}<'\cdots<'1^1 <'0<'1<'\cdots<'n.
$$
The order $<'$ is the ``good'' one to give a combinatorial interpretation of the length function in $G(r,n)$ in terms of inversions (see e.g. \cite{Re}, \cite{Ba}), while the order $<$ is often used in the study of some algebraic aspects such as the invariant theory of $B_n$ and $G(r,n)$ (see \cite{AR}, \cite{Ca1}). Here we give a motivation of these choices.
So we denote
\begin{itemize}
\item $\Des'_G(g)=\{i\in[0,n-1] \mid g(i)>'g(i+1)\}$;
\item $\des'_G(g)= |\Des'_G(g)|$;
\item $\fmaj'(g)=r\cdot \sum_{i \in \Des'_G(g)}i+\col(g)$.
\end{itemize}

If one considers the special case of Theorem \ref{desfmajcol} with $p=s=1$, and relate it to \cite[Theorem 5.1]{BZ} one can deduce that the two polynomials
$$
\sum_{g\in G(r,n)}t^{\des_G(g)}q^{\fmaj(g)}a^{\col(g)}\textrm{ and } \sum_{g\in G(r,n)}t^{\des'_G(g)}q^{\fmaj'(g)}a^{\col(g)}
$$
are equal. This can also be easily proved bijectively.
\begin{prop}\label{invol}
There exists an explicit involution $\phi:G(r,n)\rightarrow G(r,n)$ such that
$$
\Des_G(\phi(g))=\Des'_G(g)\textrm{ and }\col(\phi(g))=\col(g),
$$
and in particular $\fmaj(\phi(g))=\fmaj'(g)$.
\end{prop}
\begin{proof}
If $g\in G(r,n)$ let $S(g)=\{g(1),\ldots,g(n)\}$. The set $S(g)$ is totally ordered by both $<$ and $<'$. We let  $\iota:S(g)\rightarrow S(g)$ be the unique involution such that $\iota(x)<\iota(y)$ if and only if $x<'y$, for all $x,y\in S(g)$.
We define $\phi(g)=[\iota(g(1)),\ldots,\iota(g(n))]$: it is clear that the map $\phi$ satisfies the conditions of the statement.
\end{proof}
It follows from Proposition \ref{invol} that the two statistics $\fmaj$ and $\fmaj'$ are equidistributed on $G(r,n)$. Nevertheless the two polynomials $\sum_{g\in G(r,n)}\chi(g)q^{\fmaj(g)}$ and $\sum_{g\in G(r,n)}\chi(g)q^{\fmaj'(g)}$, where $\chi$ is any linear character of $G(r,n)$, do not coincide in general. So, in order to obtain the results appearing in \S\ref{odcafmi} we must consider the order $<$. We should also mention that the polynomial $\sum_{g\in G(r,n)}\chi(g)q^{\fmaj'(g)}$ does not factor nicely at all in general.
\smallskip

There is a more subtle difference if one considers the multivariate distributions 
$$
\sum_{g\in G(r,n)}t_1^{\des_G(g)}t_2^{\des_G(g^{-1})}q_1^{\fmaj(g)}q_2^{\fmaj(g^{-1})}a_1^{\col(g)}a_2^{\col(g^{-1})},
$$
and the analogue computed with the order $<'$.
In fact, by comparing Theorem \ref{sixgrn} and \cite[Theorem 7.1]{BZ1} we deduce that these two polynomials coincide for $r=1,2$, and are distinct if $r>2$. The case $r=1$ being trivial we can justify this coincidence for $r=2$ with a bijective proof. 
We split the set $\Des_G(g)$  into four subsets. 
$$
\Des_G(g)=\HDes_0(g)\cup \HDes_1(g)\cup \Des_{\pm}(g)\cup D_0(g),
$$
where
\begin{itemize}
\item $\HDes_0(g)=\{i\in [n-1] \mid g(i)>g(i+1)>0\}$;
\item $\HDes_1(g)=\{i\in [n-1] \mid 0>g(i)>g(i+1)\}$;
\item $\Des_{\pm}=\{i\in [n-1] \mid g(i)>0>g(i+1)\}$;
\item $D_0(g)=\{0\}$ if $g(1)<0$ and $D_0(g)=\emptyset$ otherwise.
\end{itemize}
 If we let ${\rm NN}(g)=\{i\in [n-1] \mid g(i)<0 \textrm{ and }g(i+1)<0\}$ we have
$$
\Des'_G(g)=\HDes_0(g)\cup ({\rm NN}(g)\setminus \HDes_1(g)) \cup \Des_{\pm}(g)\cup D_0(g).
$$
By the Robinson-Schensted correspondence for $B_n$ (see \cite{SW, Ca1}) we have a bijection
$$
g\mapsto[(P_0,P_1),(Q_0,Q_1)],
$$
where $(P_0,P_1)$ and $(Q_0,Q_1)$ are bitableaux of the same shape. Given a tableau $P$ let $$\Des(P)\eqdef \{i \mid \textrm{ both $i$ and $i+1$ belong to $P$ with $i$ strictly above $i+1$}\}.$$ 
If we let $\Neg(g)=\{i\in[n] \mid c_i(g)=1\}$, in the Robinson-Schensted correspondence we have:
\begin{itemize}
\item $\Neg(g)=$ content of $Q_1$, $\Neg(g^{-1})=$ content of $P_1$.
\item $\HDes_0(g)=\Des(Q_0)$, $\HDes_1(g)=\Des (Q_1)$
\item $\HDes_0(g^{-1})=\Des(P_0)$, $\HDes_1(g^{-1})=\Des (P_1)$
\item $\Des_{\pm}(g)=\{i \mid \textrm{$i$ belongs to $Q_0$ and $i+1$ belongs to $Q_1$}\}$,
\item $\Des_{\pm}(g^{-1})=\{i \mid \textrm{$i$ belongs to $P_0$ and $i+1$ belongs to $P_1$}\}$
\item $D_0(g)=\{0\}$ if and only if $1$ belongs to $Q_1$;
\item $D_0(g^{-1})=\{0\}$ if and only if $1$ belongs to $P_1$;
\end{itemize}
Let $\varphi:B_n\rightarrow B_n$ be the bijection defined by the requirement that if
\begin{eqnarray*}
g&\mapsto&[(P_0,P_1),(Q_0,Q_1)], \quad {\rm then} \\
\varphi(g)&\mapsto&[(P_0,P_1'),(Q_0,Q_1')],
\end{eqnarray*}
where $T'$ denotes the transposed tableau.
\begin{prop}
The bijection $\varphi$ satisfies the following properties:
\begin{enumerate}
\item[$(a)$] $\Neg(g)=\Neg(\varphi(g))$;
\item[$(b)$] $\Des_G(g)=\Des'_G(\varphi(g))$;
\item[$(c)$] $\Des_G(g^{-1})=\Des'_G(\varphi(g)^{-1})$.
\end{enumerate}
\end{prop}
\begin{proof}
It follows from the previous facts with the further observation that if $g\mapsto [(P_0,P_1),(Q_0,Q_1)]$ then $\Des(Q_1')={\rm NN}(g)\setminus \Des(Q_1)$ and $\Des(P_1')={\rm NN}(g^{-1})\setminus \Des(P_1)$.
\end{proof}

\begin{exa}
If $g=[5,-2,-1,-4,6,-3,-7] \in B_7$, then
$$
g\mapsto\left[\left(\begin{array}{cc}5&6\\ & \end{array}, \begin{array}{ccc}1&3&7\\2&4&\end{array}
\right),\left(\begin{array}{cc}1&5\\ & \end{array}, \begin{array}{ccc}2&4&7\\3&6&\end{array}
\right)\right].
$$
The element $\varphi(g)$ is the defined by
$$
\varphi(g)\mapsto\left[\left(\begin{array}{cc}5&6\\ & \\ & \end{array}, \begin{array}{cc}1&2\\3&4\\7&\end{array}
\right),\left(\begin{array}{cc}1&5\\ & \\ &\end{array}, \begin{array}{cc}2&3\\4&6\\7&\end{array}
\right)\right].
$$
One then can check that $\varphi(g)=[5,-3,-7,-1,6,-4,-2]$. Hence $\Neg(g)=\Neg(\varphi(g))=\{2,3,4,6,7\}$, $\Des_G(g)=\Des'_G(\varphi(g))=\{1,2,5\}$, and $\Des_G(g^{-1})=\Des'_G(\varphi(g)^{-1})=\{0,1,3,6\}$.
\end{exa}

\end{document}